\documentclass[a4 paper]{article}
\usepackage{amsthm, a4wide}
\usepackage{amsfonts}
\usepackage{amsmath}
\usepackage{amssymb}
\usepackage[english]{babel}
\newtheorem{defi}{Definition}
\newtheorem{teo}[defi]{Theorem}
\newtheorem*{rem}{Remark}

\newtheorem{coro}[defi]{Corollary}
\newtheorem{lemma}[defi]{Lemma}
\newtheorem{ex}{Example}

\newcommand{\eps}{\varepsilon}
\newcommand{\R}{I\!\!R}
\newcommand{\ox}{\overline{X}}
\newcommand{\oa}{\overline{\alpha}}
\newcommand{\ot}{\overline{\tau}}
\begin{document}
\title{Lyapunov stabilizability of controlled diffusions via a
superoptimality principle for viscosity solutions
\thanks{This research was partially supported by M.I.U.R.,
project ``Viscosity, metric, and control theoretic methods for
nonlinear partial differential equations'', and by GNAMPA-INDAM,
project ``Partial differential equations and control theory''.}}
\author{Annalisa Cesaroni\\
Dipartimento di Matematica P.  e A.\\ Universit\`a di Padova\\
via Belzoni 7, 35131 Padova, Italy\\
acesar@math.unipd.it\\ \\ }
\date{}
\maketitle

\begin{abstract}
We prove optimality principles  for semicontinuous bounded
viscosity solutions of Ha\-milton-Jacobi-Bellman equations. In
particular we provide a representation formula for viscosity
supersolutions as value functions of suitable obstacle control
problems. This result is applied  to extend the Lyapunov direct
method for stability  to controlled Ito stochastic differential
equations. We define the appropriate concept of Lyapunov function
to study the  stochastic open loop stabilizability in probability
and the local and global asymptotic stabilizability (or asymptotic
controllability). Finally we illustrate the theory with some
examples.

\smallskip \noindent{\bf Key words}.
Controlled degenerate diffusion, Hamilton-Ja\-co\-bi-Bell\-man
inequalities, viscosity solutions, dynamic programming,
superoptimality principles, obstacle problem,  stochastic control,
stability in probability, asymptotic stability.

\smallskip \noindent{\bf AMS subject classification}. 49L25, 93E15, 93D05, 93D20.
\end{abstract}

\section{Introduction}
We consider a $N$-dimensional stochastic differential equation
\[dX_t=f(X_t)dt+\sigma(X_t)dW_t\]
where $W_t$ is a standard $M$-dimensional Brownian motion. Since
the 60th, it was developed a stochastic Lyapunov method for the
analysis of the qualitative properties of the solutions of
stochastic differential equations, in analogy to the deterministic
Lyapunov method. The main contributions  are due to  Hasminskii
(see the monograph ~\cite{has} and references therein) and Kushner
(see the monograph ~\cite{ku1} and \cite{ku2}). They introduced
the notion of stability in probability and asymptotic stability in
probability. This means that  the probability the trajectory
leaves a given neighborhood of the equilibrium is decreasing to
$0$ as the initial data is approaching the equilibrium. If,
moreover, the trajectory is asymptotically approaching the
equilibrium  with probability decreasing to $0$ as the initial
data is approaching the equilibrium the system is asymptotically
stable in probability. Finally if for every  initial data the
trajectory is asymptotically approaching the equilibrium  almost
surely,  the system is asymptotically stable in the large. The
stochastic analog of deterministic Lyapunov functions $V$ are
twice differentiable continuous functions, which are positive
definite and proper  and satisfy the infinitesimal decrease
condition
\begin{equation}
-DV(x)\cdot f(x)-trace\left[a(x) D^2V(x)\right] \geq l(x) ,
\end{equation}
with $l\geq 0$ for mere Lyapunov stability and $l>0$ for $x\ne 0$
for asymptotic stability, where $a:=\sigma\sigma^T/2$. By the
Dynkin formula, this differential inequality implies that the
stochastic process $V(X_t)$, where $X_t$ is the solution of the
stochastic differential equation  starting from $x$, is a positive
supermartingale. This fact translates, in the stochastic setting,
the requirement on the Lyapunov function to decrease along the
trajectories of the dynamical system. There is a large literature
on this kind of stochastic stability: we refer to the cited
monographs and to ~\cite{mao}, see also references therein. We
recall here also the work of Florchinger \cite{flo, flo3} and
Deng, Krsti\'c, and Williams \cite {deng} on feedback
stabilization for controlled stochastic differential equations by
the Lyapunov function method.

In this paper we extend the Lyapunov method for stochastic
differential equations essentially in two directions. First of all
we consider controlled stochastic differential equations in $\R^N$
 \[dX_t=f(X_t, \alpha_t)dt+\sigma(X_t, \alpha_t)dW_t,\]
moreover  we allow the Lyapunov functions to be merely lower
semicontinuous. The nonexistence of smooth Lyapunov functions is
well known in the deterministic case, see \cite{br} for stable
uncontrolled systems and the surveys \cite{son, br} for
asymptotically stable controlled systems.  Also in the stochastic
case, the assumption of smoothness for Lyapunov functions is not
necessary and would limit considerably the applicability of the
theory and the possibility of getting a complete
Lyapunov-characterization of the stabilizability in probability by
means of a converse theorem. Kushner proved in \cite{kus1} a
characterization of asymptotic uniform stochastic stability (for
uncontrolled systems) by means of only continuous Lyapunov
functions (here, however, the infinitesimal decrease condition is
not given with a differential inequality but in terms of the weak
generator of the process). For stability in probability,
Hasminskii provided a $\mathcal{C}^2$ Lyapunov function under the
assumption of strict nondegeneracy of the diffusion: this result
cannot be extended to possibly non-degenerate diffusions. Converse
theorems in the controlled case will appear in the Ph.D. thesis by
the author \cite{ces3}. In particular we prove that the existence
of a local Lyapunov function is also necessary for the stability
in probability. Hence we show that if the system (CSDE) is
uniformly asymptotically stabilizable in probability then there
exists a local strict Lyapunov function, which is {\em
continuous}.

We define then a Lyapunov function for the stability in
probability as a {\em lower semicontinuous}, positive definite,
proper function $V$, continuous at 0 and satisfying in viscosity
sense the differential Hamilton-Jacobi-Bellman inequality
 \begin{equation}\label{eqintr}
\max_{\alpha\in A} \left\{-DV(x)\cdot
f(x,\alpha)-trace\left[a(x,\alpha) D^2V(x)\right]
\right\} \geq l(x), 
\end{equation} and we call it strict Lyapunov function if $l>0$
off 0. Our main results are the natural extensions to the
controlled diffusions of the First and Second Lyapunov Theorems:

\noindent{\em the existence of a local Lyapunov function implies
the  (open loop) stabilizability in probability of (CSDE); a
strict Lyapunov function implies the (open loop) asymptotic
stabilizability in probability}.

This means that if there exists a Lyapunov function,  then for
every initial data we can find an admissible control which keeps
the whole trajectory near the equilibrium with probability
decreasing to zero as the starting point of the trajectory
approaches the equilibrium. Moreover, the existence of a strict
Lyapunov function implies also that, for every $\eps>0$, there
exists an admissible control driving the trajectory asymptotically
to the equilibrium with probability greater than $1-\eps$. The
same proof provides the global versions as well: if $V$ satisfies
(\ref{eqintr}) in $\R^{N}\setminus\{0\}$ then $(CSDE)$ is also
(open loop) {\em Lagrange stabilizable}, i.e. has the property of
uniform boundedness of trajectories, and if $V$ is strict then the
system is  (open loop) asymptotically stabilizable in the large.
We also give sufficient conditions for the stability of viable
(controlled invariant) sets more general than an equilibrium
point.

The main tool to provide such a result is a superoptimality
principle for lower semicontinuous bounded viscosity
supersolutions $V$ of the Hamilton-Jacobi-Bellman equation
(\ref{eqintr}). The new point is that this superoptimality
principle holds as an equality and then gives a representation of
such  $V$ as  value functions of an appropriate obstacle control
problem. A similar approach has been exploited in the
deterministic case by Barron and Jensen (see \cite{baje}) for
globally asymptotically stable systems affected by disturbances
and by Soravia (\cite{sor1}, \cite{sor4}) for stable systems with
competitive controls.

Precisely we prove that

\noindent {\em every bounded LSC viscosity supersolution $V$ of
(\ref{eqintr}) can be represented as
\begin{equation}\label{fint}
V(x)= \inf_{\alpha}\sup_{t\geq 0}\mathbf{E}_x\left[
V(X_{t}^{\alpha})+\int_0^{t}\!\!l(X_s^\alpha)ds\right].
\end{equation}}

This representation formula is important on its own, since it
refers to Hamilton-Jacobi-Bellman equations for which it is not
expected uniqueness of solutions. In particular this formula
permits to give a characterization of the minimal nonnegative LSC
viscosity supersolution of the equation (\ref{eqintr}) as the
value function
\[V(x)=\inf_\alpha\mathbf{E}_x\int_0^{+\infty}\!\!l(X_s^\alpha)ds.\]

The representation formula (\ref{fint}) is obtained by introducing
a suitable sequences of obstacle problems, solved in the viscosity
sense by $V$. The conclusion comes from a uniqueness results for
viscosity solutions of such problems and from an approximation
procedure. The delicate point is the proof of a suboptimality
principle for the min-max
 value function (on the right of equality
(\ref{fint})). To get this principle, we have to choose
appropriately the  controls we are allowing for our problem.
Actually the class of controls on which we minimize a certain
given functional has to satisfy two key properties in order  to
get a dynamic programming principle. They are the stability under
concatenation and stability under measurable selection. For the
definition of the classes of controls and for the compactness and
measurable selection results we are going to use, we refer mainly
to the article by Haussmann and Lepeltier ~\cite{hl} (see also the
article by El Karoui and others ~\cite{elk} and the book by
Stroock and Varadhan ~\cite[ch 12]{sv}). For related results on
the existence of optimal controls for stochastic problems we refer
to  the article by Kushner \cite{kuscon}.

There is a large literature on dynamic programming and
superoptimality and suboptimality principles for viscosity
solutions of second order Hamilton-Jacobi-Bellman equations,
starting from the papers by P.L.Lions \cite{pll} (see also the
book \cite{fs}). We recall here the recent work by Soner and Touzi
on dynamic programming for stochastic target problems (\cite{st},
\cite{st1}). We refer also to the paper by Swiech \cite{sw} on
superoptimality and suboptimality principles for value functions
of stochastic differential games (see also the paper by Fleming
and Souganidis \cite{fs}). In this paper we are extending to the
stochastic case some results obtained by Soravia in ~\cite{sor1}
(see also ~\cite{sor2} and ~\cite{bcd}). He provides, in the
general context of differential games, a representation formula
for supersolutions of first order Isaacs equations. This gives a
superoptimality principle which holds as an equality  and which
refers to equations which in general have not unique solutions.
This result is then applied to the Lyapunov characterization of
the stabilizability to an equilibrium set of a deterministic
dynamical system with competitive controls by means of Lyapunov
functions which are only required to be continuous on the boundary
of the equilibrium set.

In the last section, we present a simple application of our
Lyapunov method. We consider an asymptotically controllable
deterministic system and we study under which conditions it
remains stable if we add to it a stochastic perturbation.  By the
Lyapunov characterization of asymptotic controllability provided
by Clarke, Ledyaev, Rifford, Stern (\cite{clrs}) and  Rifford
(\cite{r}), we know that the unperturbed system admits a Lyapunov
function, which is Lipschitz continuous and semiconcave except
possibly at the origin. We study therefore under which
perturbations this function remains a Lyapunov function also for
the perturbed system. In particular we get a small intensity
condition on the diffusion matrix $\sigma$, depending on the
semiconcavity constant of the Lyapunov function and on qualitative
properties of the stable trajectories of the deterministic
systems.

We conclude with some additional references. We recall that there
are other notions of stochastic stability. Kozin introduced the
exponential almost sure stability of uncontrolled stochastic
system. The stability in mean square and the $p$-stability were
studied by means of Lyapunov functions (we refer to the monograph
\cite{has}). In the controlled case, in previous papers Bardi and
the author (see \cite{bc1}) characterized by means of appropriate
Lyapunov functions the almost sure stabilizability of stochastic
differential equations. This is a stronger notion of stochastic
stability, never verified for nondegenerate processes. Indeed a
system is almost surely stabilizable if it behaves as a
deterministic stabilizable system and  remains almost surely in a
neighborhood of the equilibrium point. Turning to deterministic
controlled systems,  a complete Lyapunov characterization of the
asymptotic stabilizability (called asymptotic controllability) has
been proved by Sontag and Sussmann (see the articles \cite{son0},
\cite{sosu} and the review paper \cite{son}). The infinitesimal
decrease condition of the Lyapunov function along the trajectories
of the system is expressed in terms of Dini directional
derivatives, contingent directional derivatives and proximal
subgradients. There is a large literature on the stabilization of
deterministic controlled system by the Lyapunov function method:
we refer to the monograph \cite{br}, to the papers \cite{clss},
\cite{r}, see also the references therein.

The paper is organized as follows. In Section 2 we introduce the
stochastic control problems and recall the definitions and the
basic properties of the controls we are using. Section 3 is
devoted to the proof of the representation formula (\ref{fint})
for bounded continuous viscosity solutions of the differential
inequality (\ref{eqintr}) and then to the extension to lower
semicontinuous functions. Section 4 contains the definitions of
stabilizability in probability, asymptotic stabilizability and
Lyapunov functions; in Section 5 we apply the results obtained in
Section 3 to prove local and global versions of the Lyapunov
theorems.  We show that, given a control Lyapunov function or a
strict control Lyapunov function, the system (CSDE) is
respectively stabilizable or asymptotically stabilizable in
probability. In Section 6 we introduce the notion of controlled
attractor and we discuss the generalization of the direct Lyapunov
method to the case of stabilization of viable sets. Finally  in
Section 7 we present some examples illustrating the theory.

\section{Stochastic control setting}\label{sezionecontrolli}
In this section we introduce the stochastic control problem and
recall the definitions and the basic properties of the controls we
are using.

We consider a controlled  Ito stochastic differential equation:
\[(CSDE)\left\{\begin{array}{l}
dX_t=f(X_t,\alpha_t)dt +\sigma(X_t,\alpha_t)dB_t,\;\; t>0,\\
X_0=x.
\end{array}
\right. \]
We assume that $\alpha_t$ takes values in a given compact set
$A\subseteq \R^M$, $f, \sigma$ are continuous functions defined in
$\R^N\times A$, taking values, respectively, in $\R^N$ and in the
space of $N\times M$ matrices, and satisfying for all $x,y\in\R^N$
and all $\alpha\in A$
\begin{equation}\label{condition2}
\mbox{$|f(x,\alpha)-f(y,\alpha)|+
\Vert\sigma(x,\alpha)-\sigma(y,\alpha)\Vert\leq C|x-y|$. }
\end{equation}

\noindent We define
$$a(x,\alpha):=\frac{1}{2}\sigma(x,\alpha)\sigma(x,\alpha)^T$$
and assume \begin{equation}\label{convex3}
\left\{(a(x,\alpha),f(x,\alpha)) \; : \; \alpha\in A \right\}\quad
\text{is convex for all } x\in \R^N .
\end{equation}
We recall here the definition of admissible controls  that we are
allowing for our control problems. For precise definitions we
refer to \cite{hl} and \cite{elk} (see also references therein).
Actually in these articles the problem is formulated in terms of
solutions of  the martingale problem, but it is also showed that
there is an equivalent formulation in terms of solutions of
(CSDE).

We are relaxing the control problem  by using {\em weak controls},
that is, admitting all the weak solutions of (CSDE). We have not
assigned a priori a probability space
$(\Omega,\mathcal{F},\mathcal{F}_t,\mathbf{P})$ with  its
filtration. So when we introduce a control we mean that we are at
the same time choosing also a probability space and a standard
Brownian motion $B_t$ on this space. Actually under the hypothesis
(\ref{condition2}) it can be shown that the space of {\em strong
controls} is not empty and that, under suitable assumptions on the
cost functional (which are essentially the lower semicontinuity
with respect to the $x$ variable), the strong problem and the weak
problem have the same value (\cite[Theorem 4.11]{elk}).

\begin{defi}[Strict controls, Definition 2.2
\cite{hl}]\label{defsc} For every initial data $x\in\R^N$, a {\em
strict control} is a progressively measurable $A$-valued process
$(\alpha_t)_{t\geq 0}$ such that there exists a $\R^N$ valued,
right continuous, almost surely continuous, progressively
measurable solution $X_t^\alpha$ to (CSDE) (see also
~\cite[Definition 1.4]{elk}). We denote with $\mathcal{A}_x$ the
set of strict controls for $x\in\R^N$.\end{defi}

The class of strict controls can be embedded, as in the
deterministic case, in a larger class of admissible controls.  We
denote by $M(A)$ the set of probability measures on $A$ endowed
with the topology of weak convergence. We note that  it is a
separable metric space.

\begin{defi}[Relaxed controls, Definition 3.2 \cite{hl}]\label{defrc}
For every initial data $x\in\R^N$, a {\em relaxed control} is a
progressively measurable $M(A)$-valued process $(\mu_t)_{t\geq 0}$
such that there exists a $\R^N$ valued, right continuous, almost
surely continuous, progressively measurable solution $X_t^\mu$ to
(CSDE) (see also ~\cite[Definition 2.4]{elk}). We denote with
$\mathcal{M}_x$ the set of relaxed controls for $x\in\R^N$.
\end{defi}

We choose now a {\em canonical} probability space for our control
problem. By means of this canonical space we can give a
formulation of the optimization problem in a convex compact
setting. The most natural canonical space for the strict control
problem seems to be the space of the trajectories $X_{.}$ of
(CSDE). It is the space of continuous functions $\mathcal{C}$ from
$[0,+\infty)$ to $\R^N$ with its natural filtration. In order to
give a control on this space, it is sufficient to specify the
probability measure on $\mathcal{C}$ (which is the law of the
process $X_.$) and the progressively measurable function $\alpha$.
Rather than working with this canonical space we consider the
space of trajectories $(X_.,\mu_.)$ for $\mu$ relaxed control. Let
$\mathcal{V}$ the space of measurable functions from $[0,+\infty)$
to $M(A)$ with its canonical filtration. We denote with
$M(\mathcal{V})$ the set of probability measures on $\mathcal{V}$
endowed with the stable topology (this is a topology introduced by
Jacod and Menin, for precise definition we refer \cite[section
3.10]{hl} and references therein). The canonical space for the
relaxed control problem will be the product space
$\mathcal{C}\times\mathcal{V}$ with the product filtration. We
call  {\em canonic relaxed control} or {\em control rule} a
relaxed control defined in this canonical space (see
\cite[Definition 3.12]{hl} and also ~\cite[Definition 3.2]{elk}).
In order to identify  a canonic relaxed control, it will be
sufficient to specify the probability measure on the space
$\mathcal{C}\times\mathcal{V}$: the canonic relaxed controls can
be considered as measures on the canonical space.

In the following we will consider a cost functional
\[J(x,\alpha)=\sup_{t\geq 0}\mathbf{E}_x \left[V(X_t)+\int_0^t l(X_s)ds\right]\]
where $l$ is a  continuous, nonnegative function and $V$ is a LSC,
nonnegative function. The functional $j(x,t,\alpha)=\mathbf{E}_x
[V(X_t)+\int_0^t l(X_s)ds]$ satisfies, for every $x$ and $t$, the
lower continuity assumptions required in \cite{hl} on the cost
functional. Then, since the supremum of LSC maps is LSC, we get
that also the functional $J(x,\alpha)$ satisfies the same lower
semicontinuity assumptions. We list here the results obtained in
\cite{hl} that we are going to use. The crucial assumption  for
all of them, besides the right choice of the class of admissible
controls and the lower semicontinuity of the cost functional, is
the convexity assumption (\ref{convex3}).

The class of control rules is the class on which it is possible to
formulate a dynamic programming principle and to show the
existence of an optimal control. The key result is Proposition 5.2
in \cite{hl}:

\noindent {\em for every initial data $x$, the class of optimal
control rules admissible for $x$ is convex and compact.}

We have the following theorem stating the existence of an optimal
control.
\begin{teo}[Theorem 4.7 and Corollary 4.8 \cite{hl}] \label{optimalcontrol}
Under the convexity assumption (\ref{convex3}) and the other
assumptions listed above, for every initial data $x\in\R^N $ there
exists an optimal control rule for the control problem
\[\inf_{\alpha}J(x,\alpha).\]
Moreover the infimum of the cost functional computed on the class
of control rules coincides with the infimum of the cost functional
computed on the class of strict controls:
\begin{equation}\label{controllistretti}
\inf_{\alpha\in\mathcal{A}_x}J(x,\alpha)=\inf_{\alpha\in\mathcal{M}_x}J(x,\alpha).\end{equation}
In particular the optimal control can be chosen  strict.\end{teo}

The two crucial properties on the control space to get a dynamic
programming principle are  the {\em stability under measurable
selection} and the {\em stability under concatenation} (see
\cite{st1}). They are satisfied by the class of control rules. We
consider a measurable set valued map from $\R^N$ to the space of
probability measures on the canonical space
$\mathcal{C}\times\mathcal{V}$, with convex compact values. Then,
by a standard measurable selection theorem (see \cite[Theorem
5.3]{elk}) this map has a measurable selector. In ~\cite[Lemma
5.5]{hl} (see also ~\cite[ch 12]{sv} and ~\cite[Theorem 6.3,
6.4]{elk}) it is proved that this measurable selector is an
admissible control rule. Moreover, in ~\cite[Lemma 5.8]{hl} (see
also ~\cite[Theorem 6.2]{elk}) it is shown that if we take an
admissible control and then at some later stopping time we switch
to an $\eps$-optimal control from then on, the concatenated object
is still admissible.

Finally we observe that all these results remain valid if we
consider instead of the trajectories of (CSDE) in $\R^N$, the
trajectories of this system stopped at the exit time from a given
open set (see \cite{hl}).
\section{Superoptimality principles}\label{sezionerappresentazione}

In this section we prove a representation formula for bounded LSC
viscosity supersolutions of Hamilton-Jacobi-Bellman equations. We
start proving the result for continuous functions and then, by a
standard approximation procedure, we extend it to LSC bounded
functions. The representation formula is given first in the whole
space $\R^N$ and then also in a localized version. We are adapting
to the second order case the proof of optimality principles for
viscosity supersolutions of first order Hamilton-Jacobi equations
given by Soravia in ~\cite{sor1} and ~\cite{sor2}. This requires
the use of stochastic control instead of deterministic control.

We consider the following Hamilton-Jacobi-Bellman equation
\begin{equation}\label{hjb}\max_{a\in A}\left\{-f(x,a)\cdot DV(x)- trace
[a(x,a)D^2V(x)]\right\}-l(x)= 0,\end{equation} where
$l:\R^N\rightarrow\R$ is a nonnegative bounded continuous
function.

\begin{teo}[Representation formula for viscosity supersolutions]
\label{superoptimal} Consider a bounded LSC function
$V:\R^N\rightarrow \R$. If $V$ is a viscosity supersolution of the
Hamilton-Jacobi-Bellman equation (\ref{hjb}) in $\R^N$, then it
can be represented as
\begin{equation}\label{supermartingale} V(x)=
\inf_{\alpha}\sup_{t\geq
0}\mathbf{E}_x\left[V(X_{t}^{\alpha})+\int_0^{t}\!\!l(X_s^\alpha)ds\right],
\end{equation}
where the infimum is taken over all strict admissible controls.

\end{teo}
\begin{proof}
Without loss of generality, we can reduce to the case $V\geq 0$ by
an appropriate translation. Since  $V$ is LSC, bounded and
nonnegative, we can consider an increasing sequence of continuous,
nonnegative, bounded functions $V_k$ such that
\[V(x)=\sup_{k\geq 0} V_k(x)\quad\text{ for every }x\in\R^N.\] If $V$ is
continuous, we choose $V_k=V$ for every $k$.

Now for every $k\geq 0$, we introduce the following obstacle
problem in $\R^N$ with unknown $W$ and obstacle $V_k$:
\begin{align}\label{opl}
  \min \{&\lambda W(x)+ \max_{a\in A}\left[-f(x,a)DW(x)-trace
\ a(x,a) D^2W(x)\right]-l(x),\nonumber \\ & W(x)-V_k(x)\}=0.
\end{align}
Obviously $V$ is a bounded LSC viscosity supersolution of the
problem (\ref{opl}) for every $\lambda\geq 0$ and every $k\geq 0$.

For $\lambda>0$ fixed, define
\[L_{\lambda,k}(x)=\inf_{\alpha} \sup_{t\geq 0}
\mathbf{E}_x \left[e^{-\lambda t}V_k(X_{t}^{\alpha})+\int_0^{t}
 l(X_s^{\alpha})e^{-\lambda s}ds
\right].  \]

The plan of the proof is the following. First of all we show that
the value function $L_{\lambda,k}$ is a bounded discontinuous
viscosity subsolution of the obstacle problem (\ref{opl}): this
means that the upper semicontinuous envelope $L_{\lambda,k}^\star$
is a viscosity subsolution of (\ref{opl}). Then, by the comparison
principle for bounded discontinuous viscosity solutions of Isaacs
equations we get that, for every $\lambda>0$,
$L_{\lambda,k}(x)\leq V(x)$. From this, we can conclude, sending
$\lambda$ to $0$ and $k$ to $+\infty$, that $V$ satisfies the
superoptimality principle (\ref{supermartingale}).

By the definition and the boundness of $V_k$, we get that
$L_{\lambda,k}$ is bounded. We want to prove that its upper
semicontinuous envelope $L_{\lambda,k}^{\star}$ is a viscosity
subsolution of the obstacle problem (\ref{opl}). To get this
result  it is sufficient to check that $L_{\lambda, k}^\star$ is a
viscosity subsolution of the Hamilton-Jacobi-Bellman equation
\begin{equation}\label{hjb2}
\lambda L(x)+\max_{a\in A}\left\{-f(x,a)\cdot DL(x)- trace
[a(x,a)D^2L(x)]\right\}-l(x)= 0 \end{equation} at the points $x$
where $L_{\lambda, k}^\star(x)>V_k(x)$.  This result can be
obtained by standard methods in the theory of viscosity solutions
if we prove a local suboptimality principle for $L_{\lambda,
k}^\star$ on such points $x$ (see \cite{cil} and \cite{fs}).

First of all we need the following technical lemma whose proof  we
postpone to the end.
\begin{lemma}\label{soravia}
If $L^{\star}_{\lambda,k}(x)>V_k(x)$ then there exists a sequence
$x_n\rightarrow x$ with  $L_{\lambda,k}(x_n) \rightarrow
L^{\star}_{\lambda,k}(x)$ and $L_{\lambda,k}(x_n)>V_k(X_n)$ for
which there exists $\eps>0$ such that
\begin{equation}\label{lemmasoravia}
L_{\lambda,k}(x_n)\leq\inf_{\alpha} \mathbf{E}_{x_n}
\left[e^{-\lambda t }L_{\lambda,k}^{\star}(X_n ^{\alpha}(t))
+\int_0^{t} l(X_n^{\alpha}(s))e^{-\lambda s}ds
\right]\end{equation} for $t\leq \eps$ and $|x_n-x|\leq \eps$.
\end{lemma}
This is a local suboptimality principle.  This inequality, by a
standard argument in the theory of viscosity solution (for the
detailed argument see for example \cite{bar}, see also \cite{fs}),
 implies that $L^{\star}_{\lambda,k}$ is a viscosity subsolution
of equation (\ref{hjb}) at the points $x$ such that
$L^{\star}_{\lambda,k}(x)>V_k(x)$. From this we deduce that
$L_{\lambda,k }^\star$ is a viscosity subsolution of (\ref{opl}).



By the comparison principle obtained for Isaacs operators and
bounded discontinuous viscosity solutions  by Ishii in
~\cite[Theorem 7.3]{is},  we get that $V(x)\geq
L_{\lambda,k}^{\star}(x)$ for every $\lambda>0$ and $k\geq 0$. In
particular, for $T>0$ and $k$ fixed, we get
\[V(x)\geq \lim_{\lambda\rightarrow 0}L_{\lambda,k}(x)\geq
\lim_{\lambda\rightarrow 0}\inf_{\alpha} \sup_{ t\in [0,T]}
\mathbf{E}_x \left[e^{-\lambda t}V_k (X_t^{\alpha})+\int_0^t
l(X_s^{\alpha})e^{-\lambda s}ds \right]\geq \]\[\geq
\lim_{\lambda\rightarrow 0}e^{-\lambda
T}\inf_{\alpha}\sup_{t\in[0,T]}\mathbf{E}_x\left[
V_k(X_t^{\alpha})+\int_0^t\!\!l(X_s^\alpha)ds\right].\] Therefore
for every $T>0$ and $k\geq 0$
\[ V(x)\geq \inf_{\alpha}\sup_{t\in[0,T]}\mathbf{E}_x\left[
V_k(X_t^{\alpha})+\int_0^t\!\!l(X_s^\alpha)ds\right].\] Now we
want to pass to the limit for $k\rightarrow+\infty$. For $x$ fixed
and every $k\geq 0$ we consider an admissible control $\alpha_k$
such that
\begin{equation}\label{lambdak}V(x)+\frac{1}{k}\geq
\sup_{t\in[0,T]}\mathbf{E}_x\left[
V_k(X_t^{\alpha_k})+\int_0^t\!\!l(X_s^{\alpha_k})ds\right].\end{equation}
By the definitions recalled in Section 2 we can associate to each
couple $(X_.^{\alpha_k},\alpha_k)$ a control rule $P_k$. By the
compactness of the space of control rules, we can extract a
subsequence of control rules, which we continue to denote with
$P_k$, to some control rule $P$ in the stable topology. This
control rule is the measure of a trajectory of (CSDE) $(X_.,\mu)$
driven by a relaxed control $\mu$. Since the convergence in the
stable topology implies in particular the weak convergence of the
measures $P_k$ to the measure $P$, we get immediately that for
every $t\in[0,T]$
\[\lim_{k\rightarrow
+\infty}\mathbf{E}_x\int_0^t\!\!l(X_s^{\alpha_k})ds=
\mathbf{E}_x\int_0^t\!\!l(X_s) ds\] where the expected value on
the left and right hand side is computed, respectively, using the
measures $P_k$ and $P$.

Recalling now  that $V=\sup_k V_k$ where $V_k$ are continuous
functions, it is easy to show that $V$ can be obtained as
\[V(x)=\liminf_{k\rightarrow +\infty, y\rightarrow x}V_k(y):=
\sup_{\delta}\inf\left\{V_k(y)\ \left|\ |x-y|\leq\delta\ \
k\geq\frac{1}{k}\right.\right\}.\]Moreover, since the convergence
in the stable topology implies also the convergence in probability
of $X_t^{\alpha_k}$ to $X_t$, we get that, for $t\in [0,T]$ fixed,
we can extract a subsequence $X_t^{\alpha_k}$ which converges to
$X_t$ almost surely with respect to the measure $P$. Then along
this subsequence
\[\liminf_{k\rightarrow +\infty}V_k(X_t^{\alpha_k})\geq
\liminf_{k\rightarrow +\infty, y\rightarrow X_t}V_k(y)\geq
V(X_t)\quad P\text{ almost surely.}\] By the Fatou lemma and the
definition of stable convergence  we deduce that, for each $t\in
[0,T]$, along a subsequence
\[\liminf_{k\rightarrow
+\infty}\mathbf{E}V_k(X_t^{\alpha_k})\geq \mathbf{E}V(X_t),\]
where the expected value is computed respectively using the
measures $P_k$ and $P$.

To summarize, for every $t\in [0,T]$ we get from (\ref{lambdak}):
\[V(x)\geq \lim_{k\rightarrow +\infty} \mathbf{E}_x\left[
V_k(X_t^{\alpha_k})+\int_0^t\!\!l(X_s^{\alpha_k})ds\right]\geq
\mathbf{E}_x\left[ V(X_t)+\int_0^t\!\!l(X_s)ds\right].\] So for
every $T>0$ there exists a control rule for which \[V(x)\geq
\sup_{t\in[0,T]}\mathbf{E}_x\left[
V(X_t)+\int_0^t\!\!l(X_s)ds\right].\] Now, by the statement
(\ref{controllistretti}) in Theorem \ref{optimalcontrol}, we
obtain
\begin{equation}\label{lambda} V(x)\geq \inf_{\alpha}\sup_{t\in[0,T]}\mathbf{E}_x\left[
V(X_t^{\alpha})+\int_0^t\!\!l(X_s^\alpha)ds\right].\end{equation}
Now it remains only to let $T\rightarrow +\infty$.  For $\eps>0$,
consider an $\eps$ optimal control $\alpha$ for (\ref{lambda}): in
particular it gives $V(x)+\frac{\eps}{2}\geq \mathbf{E}_x\left[
V(X_T^{\alpha})+\int_0^T\!\!l(X_s^\alpha)ds\right]$. Considering
$X_T^{\alpha}$ as starting point of the trajectory we obtain by
(\ref{lambda}) \[ V(X_T^{\alpha}) \geq
\inf_{\beta}\sup_{t\in[0,T]}\mathbf{E}_{X_T^{\alpha}}\left[
V(X_{t+T}^{\beta})+\int_0^t\!\!l(X_{s+T}^\beta)ds\right]\ \ \ a.s.
\] Let $\beta$ be an $\frac{\eps}{2^2}$ optimal control rule for
$V(X_T^{\alpha})$ (we can choose a measurable selection of the
optimal  controls rules for $V(X_T^{\alpha})(\omega)$), moreover
the control rule obtained concatenating this selected control and
$\alpha$ is still an admissible control rule. We can then proceed
recursively and conclude by induction that we can construct an
admissible control rule $P$ such that
\[V(x)+\eps\geq \sup_{t\geq 0}\mathbf{E}_x\left[
V(X_t)+\int_0^tl(X_s)ds\right].\] By  the statement
(\ref{controllistretti}) in Theorem \ref{optimalcontrol} and
recalling that $\eps$ is arbitrary, we obtain
\begin{equation} \label{comp} V(x)\geq \inf_\alpha \sup_{t\geq
0}\mathbf{E}_x\left[ V(X_t)+\int_0^t\!\!l(X_s)ds\right] \geq
V(x)\end{equation} which is the desired formula. \end{proof}

We give here the proof of the technical Lemma \ref{soravia}
\begin{proof}{\it of Lemma \ref{soravia}.}
If the statement were not true, for every  sequence
$x_n\rightarrow x$ such that $L_{\lambda,k}(x_n) \rightarrow
L^{\star}_{\lambda,k}(x)$ and $L_{\lambda, k}(x_n)>V_k(x_n)$ for
every $n$ we could find $t_n\leq \frac{1}{n}$ such that for
$|x_n-x|\leq \frac{1}{n}$
\begin{equation} \label{claimlemma}
L_{\lambda,k}(x_n)>\inf_{\beta} \mathbf{E}_{x_n} \left[e^{-\lambda
t_n}L_{\lambda, k}^{\star} (X_{n}^{\beta}(t_n)) +\int_0^{t_n}
l(X_n^{\beta}(s))e^{-\lambda s}ds \right].\end{equation}
 By definition  of
$L_{\lambda, k}$,
for every $\eps_n>0$  and every control $\alpha$, there exists
$T(\eps_n,\alpha)$ such that \begin{equation}\label{ineq}
L_{\lambda, k}(x_n)-\eps_n< \mathbf{E}_{x_n} \left[e^{-\lambda
T(\eps_n,\alpha)}V_k(X_n^{\alpha}(T(\eps_n,\alpha)))
+\int_0^{T(\eps_n,\alpha)} l(X_n^{\alpha}(s))e^{-\lambda s}ds
\right]. \end{equation} By  the inequality (\ref{claimlemma}), we
can choose a sequence $\eps_n\rightarrow 0$ and controls $\beta_n$
for which \begin{equation}\label{ineq2} L_{\lambda,
k}(x_n)-2\eps_n
\geq\mathbf {E}_{x_n} \left[e^{-\lambda t_n} L_{\lambda,
k}^{\star}(X_n^{\beta_n}(t_n))+\int_{0}^{t_n} l(X_n^{\beta_n}(s))
e^{-\lambda s}ds\right].\end{equation} Therefore for every control
$\alpha$ we obtain from the inequalities (\ref{ineq}) and
(\ref{ineq2})
\[  \mathbf{E}_{x_n}
\left[e^{-\lambda t_n} L_{\lambda,
k}^{\star}(X_n^{\beta_n}(t_n))+\int_{0}^{t_n} l(X_n^{\beta_n}(s))
e^{-\lambda s}ds\right]+\eps_n\leq L_{\lambda,
k}(x_n)-\eps_n<\]\[< \mathbf{E}_{x_n} \left[e^{-\lambda
T(\eps_n,\alpha)}V_k(X_n^{\alpha}(T(\eps_n,\alpha)))
+\int_0^{T(\eps_n,\alpha)} l(X_n^{\alpha}(s))e^{-\lambda s}ds
\right].
\]
We claim now that for every $n$ there exists $\alpha$ such that
$T(\eps_n,\alpha)\leq t_n$. Assume by contradiction that there
exists $N$ such that,
 for every $\alpha$ admissible, $T(\eps_N,\alpha)>t_N$,
in particular $T(\eps_N,\alpha_N)>t_N$ for every  control
$\alpha_N$ which for $t\leq t_N$  coincides with $\beta_N$: by the
previous inequality we get
\[\mathbf{E}_{x_N}
L_{\lambda, k}^{\star}(X_{N}^{\beta_N}(t_N))+\eps_N<
\mathbf{E}_{x_N} \left[e^{-\lambda
(T(\eps_N,\alpha_N)-t_N)}V_k(X_N^{\alpha_N}
(T(\eps_N,\alpha_N)))+\right.\]\[\left.
+\int_{t_N}^{T(\eps_N,\alpha_N)} l(X_N^{\alpha_N}(s))e^{-\lambda
(s-t_N)}ds \right].
\]
We choose now an admissible control rule  which is $\eps_N/2$
optimal for $L_{\lambda, k}^{\star}(X_{N}^{\beta_N}(t_N))$ (we
operate a measurable selection between the $\eps_N/2$ optimal
control rules  for $L_{\lambda,
k}^{\star}(X_{N}^{\beta_N}(t_N))(\omega)$)  and the concatenated
control with $\beta_N$ is still an admissible control rule $P_N$
which is the measure associated  to the couple $(X_N(.),\mu_N)$.
Therefore we get
\[\sup_{t\geq t_N}\mathbf{E}_N
\left[e^{-\lambda t}V_k(X_N(t)) +\int_{t_N}^{t}\!\!
l(X_N(s))e^{-\lambda s}ds \right]+\eps_N/2\leq
\]\[\leq
\mathbf{E}_N \left[e^{-\lambda T(\eps_N,\mu_N)}
V_k(X_N(T(\eps_N,P_N)))+\int_{t_N}^{T(\eps_N,\mu_N)}
\!\!l(X_N(s))e^{-\lambda s}ds \right]\] where the expected value
is computed with respect to the measure $P_N$, and the inequality
we obtain gives a contradiction. Therefore there exists for every
$n$ an admissible
 control rule $P_n$ such that $T(\eps_n,\mu_n)\leq t_n$:
choosing $\alpha=\mu_n$ in the inequality (\ref{ineq}) we get
\[L_{\lambda, k}(x_n)-\eps_n\leq \mathbf{E}_{x_n}
\left[e^{-\lambda T(\eps_n,\mu_n)}V_k(X_n(T(\eps_n,\mu_n)))
+\int_0^{T(\eps_n,\mu_n)}\!\!\ l(X_n(s))e^{-\lambda s}\right].\]
For every $n$ we call $A_n=\{\omega|\ \ X_n(T(\eps_n,\mu_n)) \in
B(x_n,(\frac{1}{\sqrt[4]{n}})\}$ and $B_n=\Omega\setminus A_n$:
since for every $n$ the trajectory is a Markov process and the
drift and the diffusion of this control problem are equi-Lipschitz
and equi-bounded in the compacts with respect to the control it is
possible to show (we refer to ~\cite[pp 284,285]{doob} for the
proof) that $\mathbf{P}(B_n)\leq K\mathcal{O}
(\frac{1}{\sqrt[4]{n}})^3$ where
$\mathcal{O}(\frac{1}{\sqrt[4]{n}})^3$ is uniform with respect to
the initial data $x_n$ and to the control $\mu_n$. Therefore we
get
\[L_{\lambda, k}(x_n)-\frac{\eps_n}{2}<
\int_{A_n}\!\!\left[V_k(X_n(T(\eps_n,\mu_n)))+ \int_0
^{T(\eps_n,\mu_n)}\!\! l(X_n(s))e^{-\lambda s} ds \right]+\]\[
+\int_{B_n}\!\!\left[V_k(X_n( T(\eps_n,\mu_n)))+
\int_0^{T(\eps_n,\mu_n)} l(X_n(s))e^{-\lambda s} ds \right]
\leq\]\[ \leq \left[\sup_{B(x,\frac{2}{\sqrt{n}})}V_k(y)+ \sup_{
B(x,\frac{2}{\sqrt{n}})} l(y) t_n\right]+ o(\frac{1}{\sqrt{n}})
\]
from which, since   $V_k$ is continuous and $t_n\leq 1/n$, letting
$n\rightarrow+\infty$,  we deduce
\[L_{\lambda, k}^{\star}(x)\leq V_k(x)\] in contradiction with our assumption.
\end{proof}
\begin{rem}\upshape
The previous result can be proved in more general situations:
consider a bounded, nonnegative, LSC viscosity supersolution
$V:\R^N\rightarrow \R$ of
\[\max_{a\in A}\left\{-f(x,a)\cdot DV(x)- trace
[a(x,a)D^2V(x)]\right\}+k(x)V(x)\geq l(x)\] where
$k:\R^N\rightarrow\R$ is a Lipschitz continuous nonnegative
function. The proof of Theorem \ref{superoptimal} applies directly
and we obtain the representation formula
\[V(x)=\inf_{\alpha}\sup_{t\geq 0}\mathbf{E}_{x}\left[
V(X_t^{\alpha})e^{-\int_0^t k(X_s^{\alpha})ds}+\int_0^t\!\!
l(X_s^\alpha)e^{-\int_0^s k(X_u^{\alpha})du}ds\right].\]

\end{rem}
We can prove also a localized version of the Theorem
\ref{superoptimal}.
\begin{coro}\label{cor}
Consider an open set $\mathcal{O}\subseteq \R^N$. For every
$\delta>0$, consider the set
$\mathcal{O}_\delta:=\{x\in\mathcal{O}\ |\ d(x,\partial
\mathcal{O})>\delta\}$ and denote with $\tau_\delta^{\alpha}$ the
stopping time at which the sample function of the process
$X_t^{\alpha}$ reaches the boundary $\partial\mathcal{O}_\delta$:
we denote with $\tau_\delta^{\alpha}(t)$ the minimum between
$\tau_\delta^{\alpha}$ and $t$. Assume that
$V:\overline{\mathcal{O}}\rightarrow\R$ is a  bounded nonnegative
function. If $V$ is a LSC viscosity supersolution of the
Hamilton-Jacobi-Bellman equation (\ref{hjb}) in $\mathcal{O}$,
then it can be represented, for every $\delta$,
$x\in\mathcal{O}_\delta$, as
\[
V(x)= \inf_{\alpha}\sup_{t\geq 0}\mathbf{E}_x\left[
V(X_{\tau_\delta^{\alpha}(t)}^{\alpha})+\int_0^{\tau_\delta^{\alpha}(t)}\!\!
l(X_s^\alpha)ds\right].
\]
\end{coro}
\begin{proof}
We fix $\delta>0$ and a smooth cut off  function $0\leq\xi\leq 1$
such that $\xi(x)=0$ for $x\in\R^N\setminus \mathcal{O}$ and
$\xi(x)=1$ for $x\in\mathcal{O}_\delta$. We consider the
stochastic controlled differential equation in $\R^N$:
\[(CSDE)'\left\{\begin{array}{l}
dX_t=f(X_t,\alpha_t)\xi^2(X_t)dt +\sigma(X_t,\alpha_t)\xi(X_t)dB_t,\;\; t>0,\\
X_0=x.
\end{array}
\right. \] Observe that for $x\in\mathcal{O}_\delta$, the solution
$(X')^\alpha$ to (CSDE)$'$ coincides a.s. with the solution
$X^\alpha$ to (CSDE) up to time $\tau_\delta^\alpha$. We define
the process $X_{\tau_\delta^{\alpha}(t)}^{\alpha}$ obtained by
stopping the process $(X')_t^{\alpha}$ at the instant it reaches
the boundary of $\mathcal{O}_\delta$: it has a Ito stochastic
differential and  it is still a strong Markov process with
continuous trajectories (see for example ~\cite[Lemma 3.3.1]{has}
and references therein).

We extend $V$ outside $\mathcal{O}$ as a bounded nonnegative LSC
function that we continue to denote $V$. So it is immediate to
show that $V$ is a viscosity supersolution in $\R^N$ of the
equation:
\begin{equation}\label{hjb3}\max_{a\in A}\left\{-f(x,a)\xi^2(x)\cdot DV(x)- trace
[a(x,a)\xi^2(x)D^2V(x)]\right\}-l(x)\xi^2(x)= 0.\end{equation} We
can apply Theorem \ref{superoptimal} to $V$. Indeed it is
sufficient  to define
\[L_{\lambda,k}(x)\!\!\!=\!\!\!\left\{\begin{array}{ll}\!\!\!\!\inf_{\alpha} \sup_{t\geq 0}
\mathbf{E}_x \left[e^{-\lambda
\tau_\delta^{\alpha}(t)}V_k(X_{\tau_\delta^{\alpha}(t)}^{\alpha})+\int_0^{t}
 l(X_{\tau_\delta^{\alpha}(s)}^{\alpha})\xi^2(X_{\tau_\delta^{\alpha}(s)}^{\alpha})
 e^{-\lambda \tau_\delta^{\alpha}(t)}ds\right]&\!\!\!\!\!\!\!\text{in }
\overline{\mathcal{O}_\delta}\\\!\!\!\! V_k(x)&
\!\!\!\!\!\!\!\text{in }\R^N\setminus
\overline{\mathcal{O}_\delta}.
\end{array}\right.\] We can repeat the proof in Theorem
\ref{superoptimal} (all the results in \cite{hl} hold also for the
stopped process $Y^\alpha$) and we get that $L^\star_{\lambda,k}$
is a viscosity supersolution of the obstacle problem (\ref{opl})
in $\R^N$. So again repeating the same arguments of Theorem
\ref{superoptimal} we get that $V$ satisfies the following
representation formula for $x\in\mathcal{O}_\delta$:
\[
V(x)= \inf_{\alpha}\sup_{t\geq 0}\mathbf{E}_x\left[
V(X_{\tau_\delta^{\alpha}(t)}^{\alpha})+\int_0^{\tau_\delta^{\alpha}(t)}\!\!
l(X_s^\alpha)ds\right].\]
\end{proof}

\begin{rem}[Minimal nonnegative solution]\upshape \label{minimalsupersolution}
These representation formulas for viscosity solutions are
interesting on their own, as we have pointed out in the
introduction: indeed they apply to Hamilton-Jacobi-Bellman
equations for which there are no comparison principles and then no
uniqueness of solutions.

We consider the following Hamilton-Jacobi-Bellman equation  in
$\R^N$
\[\max_{a\in A}\left\{-f(x,a)\cdot DV(x)- trace
[a(x,a)D^2V(x)]\right\}=l(x)
\] with $l\geq 0$: since the constant function $U\equiv 0$ is always a
subsolution, it is interesting to characterize the minimal
nonnegative supersolution.

From a control point of view, the natural solution seems to be the
value function of the infinite horizon control problem with
running cost $l$
\[V_{\infty}(x)=\inf_{\alpha} \mathbf{E}_x
\int_0^{+\infty}\!\!l(X_s^\alpha)ds.\] If $V_\infty$ is well
defined and bounded, then it is possible to show that it is LSC,
by an argument based on the properties of the class of admissible
relaxed controls. Moreover, by standard methods in the theory of
viscosity solutions (see \cite{fso},\cite{cil}), it is possible to
show that $V_\infty$  is a viscosity supersolution of the previous
Hamilton-Jacobi-Bellman equation.

In this case, an easy application of the previous theorems gives
that every bounded, nonnegative, viscosity supersolution $V$ of
the Hamilton-Jacobi-Bellman equation in $\R^N$ satisfies
\[V(x)\geq V_{\infty}(x) \] therefore
$V_\infty$ is the minimal nonnegative viscosity supersolution of
the equation. \end{rem}
\begin{rem}[Representation formula
for viscosity subsolutions]\upshape The counterpart of Theorem
\ref{superoptimal} for viscosity subsolutions is straightforward
from classical suboptimality principles: let $U:\R^N\rightarrow
\R$ be an upper semicontinuous bounded viscosity subsolution of
the Hamilton-Jacobi-Bellman equation
\[\max_{a\in A}\left\{-f(x,a)\cdot DU(x)- trace
[a(x,a)D^2U(x)]\right\}\leq l(x),
\] then the function $U$
can be represented as:
\[U(x)=\inf_{\alpha} \inf_{t\geq  0}   \mathbf{E}_x\left[
U(X_t^{\alpha})+ \int_0^t\!\!l(X_s^\alpha)ds\right].\]
\end{rem}

\section{Stability in probability and Lyapunov functions}
We begin this section with the notion of both Lyapunov and
asymptotic stability in probability. They were introduced  by
Hasminskii and Kushner (see ~\cite{has} and ~\cite{ku1}) in the
case of uncontrolled stochastic differential equations. We present
their natural extension to the case of controlled diffusions.
\begin{defi}[stabilizability in probability]\label{stabinprobdefi}

The controlled system (CSDE) is {\em (open loop) stabilizable in
probability} at the origin if for all $\eps,k>0$ there exists
$\delta>0$ such that for every $|x|\leq \delta$ there exists a
control $\oa_.\in\mathcal{A}_x$ such that the corresponding
trajectory $\ox_.$ verifies
\[\mathbf{P}_x\left(\sup_{t\geq 0}|\ox_t|\geq k\right)\leq \eps.
\] This is equivalent to assume that for every positive  $k$
\[\lim_{x\rightarrow 0} \inf_{\alpha}
\mathbf{P}_x\left(\sup_{t\geq 0}|X_t^\alpha|\geq k\right)=0.\]

The system is \emph{(open-loop) Lagrange stabilizable in
probability}, or it has the property of \emph{uniform boundedness
of trajectories}, if for each $\eps>0,R>0$ there is $S>0$ such
that, for any initial point $x$ with $|x|\leq R$,
\[\inf_{\alpha}\mathbf{P}_x\left( \sup_{t\geq 0}|X_t^\alpha|\geq
S\right)\leq \eps.\] This is equivalent to assume that for every
$R>0$
\[\lim_{S\rightarrow +\infty}\sup_{|x|\leq R} \inf_{\alpha}
\mathbf{P}_x\left(\sup_{t\geq 0}|X_t^\alpha|\geq S\right)=0 .\]

\end{defi}
\begin{rem}\upshape
The stabilizability in probability implies that the origin is a
{\em controlled
equilibrium} of $(CSDE)$, i.e., 
$$
\exists \,\oa\in A \,:\; f(0,\oa)=0, \; \sigma(0,\oa)=0.
$$
In fact, the definition gives for any $\eps>0$, for $k>0$ fixed,
an admissible control such that the corresponding trajectory
starting at the origin satisfies $\mathbf{P}(\sup_{t\geq 0}
|X_t|\geq k)\leq\eps $
so \[\mathbf{E}_0\int_0^{+\infty} l(|X_t|)e^{-\lambda t}dt \leq
\frac{\eps}{\lambda}\] for any $\lambda>0$ and any
 real function $l$ such that $0\leq l(r)\leq 1$ for any $r$
and $l(r)=0$ for $r\leq k$. Then
$\inf_{\alpha_.\in {\cal A}_0} \mathbf{E}_0\int_0^{+\infty}
l(|X_t|)e^{-\lambda t}dt=0 $.
Theorem \ref{optimalcontrol} implies that the $\inf$ is attained:
therefore for any $k>0$ there is a minimizing control which
produces a trajectory satisfying a.s. $|X_t|\leq k $ for all
$t\geq 0$. So
$\inf_{\alpha_.\in {\cal A}_0} \mathbf{E}_0\int_0^{+\infty}
|X_t|e^{-\lambda t}dt=0$ for any $\lambda>0$. Again Theorem
\ref{optimalcontrol} implies that the $\inf$ is attained, and the
minimizing control produces a trajectory satisfying a.s. $|X_t|=0$
for all $t\geq 0$. The conclusion follows from standard properties
of stochastic differential equations.

Regarding the Lagrange stabilizability, we  observe that, using
standard properties of diffusions under the regularity assumptions
(\ref{condition2}), it is possible to prove (see \cite{doob, has})
that for every fixed $T>0$ and $R>0$
\[\lim_{S\to +\infty}\sup_{|x|\leq R}\inf_{\alpha}\mathbf{P}_x\left(
\sup_{0\leq t\leq T}|X_t^\alpha|\geq S \right)=0.\] Nevertheless
the Lagrange stabilizability is a stronger condition since it
requires that
\[\lim_{S\to +\infty }\sup_{|x|\leq R}\inf_{\alpha}\sup_{T\geq 0}\mathbf{P}_x\left(
\sup_{0\leq t\leq T}|X_t^\alpha|\geq S \right)=0.\]
\end{rem}

The controlled diffusion is said to be asymptotically stabilizable
in probability if the equilibrium point is not only stabilizable
but also an attracting point for the system, locally around the
equilibrium point.
\begin{defi}[asymptotic stabilizability in probability]
The controlled system is locally {\em asymptotically stabilizable
in probability} at the origin if for all $\eps,k>0$ there exists
$\delta>0$ such that for every $|x|\leq \delta$ there exists a
control $\oa_.\in\mathcal{A}_x$ such that the corresponding
trajectory $\ox_.$ verifies
\[\mathbf{P}_x\left(\sup_{t\geq 0}|\ox_t|\geq k\right)\leq \eps
\ \ \mbox{ and }\ \ \mathbf{P}_x\left(\limsup_{t\rightarrow
+\infty}|\ox_t|>0\right)\leq \eps.
\]
This is equivalent to assume that  for all $k>0$
\[\lim_{x\rightarrow 0} \inf_{\alpha}
\left[\mathbf{P}_x\left(\sup_{t\geq 0}|X_t^\alpha|\geq k\right)+
\mathbf{P}_x\left(\limsup_{t\rightarrow
+\infty}|X_t^\alpha|>0\right)\right]=0.
\]
\end{defi}
There is a global version of the previous stability notion:
\begin{defi}[asymptotic stabilizability in the large]
The controlled system is {\em asymptotic stabilizable in the
large} at the origin if it is Lyapunov stabilizable in probability
around the equilibrium and  for every $x\in\R^N$
\[\inf_{\alpha}\mathbf{P}_x\left(
\limsup_{t\rightarrow +\infty}|X_t^\alpha|>0\right)=0.
\]
This means that for every $\eps>0$ and for every initial data $x$
we can choose an admissible control in $\mathcal{A}_x$ which
drives the trajectory to the equilibrium with probability greater
than $1-\eps$. \end{defi}

Next we give the appropriate definition of control Lyapunov
functions for the study of the stochastic stabilities defined
above.
\begin{defi}[Lyapunov function]\label{liap}
Let ${\cal O}\subseteq \R^N$ be a bounded open set containing the
origin. A function $V:{\cal O}\rightarrow \R$ is a \emph{local
Lyapunov function} for $(CSDE)$ if it satisfies the following
conditions:

\noindent(i) it is lower semicontinuous and continuous at the
origin;

\noindent(ii) it is {\em positive definite}, i.e., $V(0)=0$ and
$V(x)>0$ for all $x\neq 0$;

\noindent (iii) it is bounded;

\noindent(iv) it is a viscosity supersolution of the equation
\begin{equation} \label{eq1}
\max_{\alpha\in A}\left\{-DV(x)\cdot
f(x,\alpha)-trace\left[a(x,\alpha) D^2V(x)\right] \right\}=0 \quad
\mbox{in } {\cal O}.
\end{equation}
\end{defi}
We introduce the notion of strict Lyapunov function both in the
local and global setting.
\begin{defi}[local strict Lyapunov function]\label{sliap}
Let ${\cal O}\subseteq \R^N$ be a bounded open set containing the
origin. A function $V:{\cal O}\rightarrow \R$ is a \emph{local
strict Lyapunov function} for $(CSDE)$ if it satisfies the
conditions (i),(ii), (iii) in the previous definition and

\noindent(iv') it is a viscosity supersolution of the equation
\begin{equation} \label{eq2}
\max_{\alpha\in A}\left\{-DV(x)\cdot
f(x,\alpha)-trace\left[a(x,\alpha) D^2V(x)\right] \right\}=l(x)
\quad \mbox{in } {\cal O},
\end{equation} where $l:{\cal O}\rightarrow \R$ is a positive definite, bounded
and uniformly continuous function.
\end{defi}
\begin{defi}[global strict Lyapunov function]\label{gsliap}
Let $\mathcal{O}\subseteq\R^N$ be an open set containing the
origin.
 A function $V:\mathcal{O}\rightarrow \R$ is a
\emph{global strict Lyapunov function} for $(CSDE)$ if it
satisfies the following conditions:

\noindent (i) it is lower semicontinuous and continuous at the
origin;

\noindent(ii) it is {\em positive definite}, i.e., $V(0)=0$ and
$V(x)>0$ for all $x\neq 0$;

\noindent(iii) it is {\em proper}, i.e.,
$\lim_{x\rightarrow\partial\mathcal{O}}V(x)=+\infty$, or,
equivalently, its level sets $\{x| V(x)\leq\mu\}$ are bounded;

\noindent(iv) it is a viscosity supersolution of the equation
\begin{equation} \label{eq3}
\max_{\alpha\in A}\left\{-DV(x)\cdot
f(x,\alpha)-trace\left[a(x,\alpha) D^2V(x)\right] \right\}=l(x)
\quad \mbox{in } {\cal O},
\end{equation} where $l:\mathcal{O}\rightarrow \R$ is a positive definite
uniformly continuous function.
\end{defi}

\section{Direct Lyapunov theorems}\label{teoremidiretti}
In this section we develop a direct Lyapunov method for the study
of stabilizability in probability of controlled diffusions both in
the local and global setting. For the uncontrolled case, the
extension of the Lyapunov second method to the case of stochastic
systems is due to Hasminskii and Kushner independently (see the
monographs \cite{has}, \cite{ku1}, see also the references therein
for earlier related results).

The main tool of the proof of the Lyapunov theorems is the
representation formula for viscosity solutions obtained in Section
\ref{sezionerappresentazione}.
\begin{teo}[Stabilizability in probability]\label{liap1}
Assume conditions (\ref{condition2}), (\ref{convex3}) and the
existence of  a local Lyapunov function $V$ in the open set ${\cal
O}$. Then:

\noindent (i) the system  is stabilizable in probability,

\noindent(ii) if in addition the Lyapunov function is global, then
the system is also Lagrange stabilizable in probability.
\end{teo}
\begin{proof}
We start proving (i). We fix $k>0$ such that $B_k\subset {\cal
O}$. We fix $\eps>0$ and  define $\eta=\eps\  {\min_{|y|\geq
k}V(y)}$. We denote with $\tau_{k}^{\alpha}(x)$ the first exit
time of the trajectory $X_t^{\alpha}$  from the open ball $B_k$
centered at the origin with radius $k$. By the continuity at the
origin of $V$ we can find $\theta>0$ such that if $|x|\leq \theta$
then $V(x)\leq\frac{\eta}{2}$. The superoptimality principle in
Corollary  \ref{cor} gives, for $|x|\leq \theta\wedge k$,
\[\eta/2\geq V(x)=\inf_{\alpha}\sup_{t\geq 0}
\mathbf{E}_x V(X_{t\wedge \tau_k^{\alpha}(x)}).\] We choose now an
$\frac{\eta}{2}$ optimal control $\oa\in \mathcal{A}_x$ for the
previous control problem, we denote by $\ox_t$ the corresponding
trajectory, stopped at the exit time from ${\cal O}$, and we get
for every  $t\geq 0$
\[  \eta\geq
\mathbf{E}_x V(\ox_{t\wedge\tau_{k}^{\oa}})\geq
\int_{\{\sup_{0\leq s\leq t} |\ox_s|\geq
k\}}V(\ox_{\tau_{k}^{\oa}})\geq \mathbf{P}(\sup_{0\leq s\leq t}
|\ox_s|\geq k) \min_{|y|\geq k}V(y).\] As $t\rightarrow+\infty$,
we obtain the following bound on the probability that the
trajectory $\ox_t$  leaves the ball $B_k$
\[\mathbf{P}(\sup_{t\geq 0} |\ox_t|\geq k)\leq
\frac{\eta}{\min_{|y|\geq k}V(y)}=\eps.\] This proves the
stabilizability in probability.

We pass now to (ii). Repeating the argument above we get that for
every $k>0$
\[\inf_\alpha\mathbf{P}(\sup_{t\geq 0} |X^\alpha_t|\geq k)\leq
\frac{V(x)}{\min_{|y|\geq k}V(y)}.\] This implies the Lagrange
stabilizability: indeed given $R>0$ and $\eps>0$, we choose $k$
such that \[\frac{\max_{|y|\leq R}V(y)}{\min_{|y|\geq k}V(y)}\leq
\eps.\]
\end{proof}

In the case the system admits a strict  Lyapunov function we prove
that there exists a control which not only stabilizes the
diffusion in probability but  also drives it asymptotically to the
equilibrium. We obtain the result using standard martingale
inequalities; in the uncontrolled case, a similar proof of
asymptotic stability has been given in \cite{deng} (see also
\cite{mao}).

\begin{teo}[Asymptotic stabilizability]\label{liap2}
Assume conditions (\ref{convex3}), (\ref{condition2}) and  the
existence of a local strict Lyapunov function in an open set
$\mathcal{O}$. Then

\noindent (i) the system (CSDE)  is  locally asymptotically
stabilizable in probability;

\noindent (ii) if the strict Lyapunov function is global, then the
system (CSDE) is  asymptotically stabilizable in the large.
\end{teo}
\begin{proof} We start proving (i).  For every
$k>0$, such that $B_k\subset {\cal O}$, we get, by  Corollary
\ref{cor}, that the function $V$ satisfies, for $x\in B_k$, the
following superoptimality principle
\begin{equation}\label{asintotica}
V(x)= \inf_{\alpha}\sup_{t\geq 0}\mathbf{E}_x\left[
V(X_{\tau_{k}^\alpha(t)}^{\alpha})+
\int_0^{\tau_{k}^\alpha(t)}\!\!l(X_s^\alpha)ds\right]
\end{equation}
where the trajectories are stopped at the exit time from $B_k$. By
Theorem \ref{optimalcontrol} there exists an optimal control
$\oa\in\mathcal{A}_x$ for this value problem. We indicate with
$\ox_.$ the corresponding trajectory and with $\ot$ the exit time
from the open ball  $B_k$. Repeating the proof of Theorem
\ref{liap1} we get the stabilizability in probability:
\[\mathbf{P}_x(\ot<+\infty)=\mathbf{P}_x(\sup_{t\geq 0} |\ox_t|\geq k)\leq
\frac{V(x)}{\min_{|y|\geq k}V(y)}.\] We denote by $B(x)=\{\omega\
|\ |\ox_t(\omega)|\leq k \ \ \forall t\geq 0 \}$. By the previous
estimate
$\mathbf{P}_x(B(x))=\mathbf{P}_x\left(\ot=+\infty\right)\geq
1-V(x)/(\min_{|y|\geq k}V(y))$.

We claim that $l(\ox_t(\omega))\rightarrow 0$ as $t\rightarrow
+\infty$ for almost all $\omega\in B(x)$, from this, using the
positive definiteness of the function $l$, we can deduce that
\[\mathbf{P}_x\left(\limsup_{t\rightarrow
+\infty}|\ox_t|>0\right)\leq \frac{V(x)}{\min_{|y|\geq k}V(y)}\]
which gives, by the continuity at the origin of the function $V$,
the asymptotic stabilizability in probability.

We assume by contradiction that the claim is not true: then there
exists $\eps>0$, a subset $\Omega_\eps\subseteq B(x)$ with
$\mathbf{P}(\Omega_\eps)>0$, and for every $\omega\in
\Omega_\eps$ 
a sequence $t_n(\omega)\rightarrow +\infty$ such that
$l(\ox_{t_n}(\omega)) >\eps$.  We define
 $$
 F(k):=\max_{|x|\leq k, \alpha\in A}|f(x,\alpha)| ,
 \qquad\Sigma(k)=\max_{|x|\leq k, \alpha\in A} \|\sigma(x,\alpha)\|.
$$
 We indicate with $\ot(s)$ the minimum between $\ot$
 and $s$ and compute, for $t\geq 0$ fixed
\[
\mathbf{E}\left\{\sup_{t\leq s\leq
t+h}|\ox_{\ot(s)}-\ox_{\ot(t)}|^2 \right\}=
\]
\[
=\mathbf{E}\left\{\sup_{t\leq s\leq t+h}|\int_{\ot(t)}^{\ot(s)}
f(\ox_u,\oa_u)du+\int_{\ot(t)}^{\ot(s)} \sigma(\ox_u,\oa_u)dB_u|^2
\right\}\leq
\]
\[
\leq  2\mathbf{E}\left\{\sup_{t\leq s\leq
t+h}|\int_{\ot(t)}^{\ot(s)} f(\ox_u,\oa_u)du|^2\right\}+2
\mathbf{E}\left\{\sup_{t\leq s\leq t+h}|\int_{\ot(t)}^{\ot(s)}
\sigma(\ox_u,\oa_u)dB_u|^2 \right\}\leq
\]
\[
\leq 2 F^2(k)h^2+2 \mathbf{E}\left\{\sup_{t\leq s\leq t+
h}|\int_{\ot(t)}^{\ot(s)} \sigma(\ox_u,\oa_u)dB_u|^2\right\}.
\]
By Theorem 3.4 in \cite{doob} (the process $|\int_t^s
\sigma(\ox_u,\oa_u)dB_u|$ is a positive semimartingale) we get
\[
\mathbf{E}\left\{\sup_{t\leq s\leq
t+h}|\ox_{\ot(s)}-\ox_{\ot(t)}|^2 \right\} \leq 2 F^2(k)h^2+8
\sup_{t\leq s\leq t+h}\mathbf{E}\left\{|\int_{\ot(t)}^{\ot(s)}
\sigma(\ox_u,\oa_u)dB_u|^2\right\}\leq
\]
\[
 \leq 2
F^2(k)h^2+8 \sup_{t\leq s\leq
t+h}\mathbf{E}\left\{\int_{\ot(t)}^{\ot(s)}
|\sigma(\ox_u,\oa_u)|^2du\right\}\leq 2 F^2(k)h^2+8 \Sigma^2(k)h .
\]
Then, Chebyshev inequality gives, for every $t\geq 0$ fixed
\begin{equation}\label{continuitytrajectory}
\mathbf{P}\left(\sup_{t\leq s\leq t+h}
|\ox_{\ot(s)}-\ox_{\ot(t)}|>r \right) \leq  \frac{2 F^2(k)h^2+8
\Sigma^2(k)h}{r^2}.\end{equation}
Since $l$ is continuous, we can fix $\delta$ such that
$|l(x)-l(y)|\leq\frac{\eps}{2}$ if $|x-y|\leq\delta$ and $|x|, |y|
\leq k$: we compute
\[
\mathbf{P}\left(\sup_{t\leq s\leq t+h}
|l(\ox_{\ot(s)})-l(\ox_{\ot(t)})| \leq\frac{\eps}{2}\right)\geq
\mathbf{P}\left(\sup_{t\leq s\leq t+h}
|\ox_{\ot(s)}-\ox_{\ot(t)}|\leq \delta\right)\geq\]\[\geq
1-\frac{2 F^2(k)h^2+8 \Sigma^2(k)h}{\delta^2}.
\] We choose $h$ such that $0<\frac{2 F^2(k)h^2+8
\Sigma^2(k)h}{\delta^2}\leq \mathbf{P}_x(\Omega_\eps)-r$ for some
$r>0$ so that for every $t\geq 0$
\begin{equation}\label{disu1}
\mathbf{P}\left(\sup_{t\leq s\leq t+h}
|l(\ox_{\ot(s)})-l(\ox_{\ot(t)})| \leq\frac{\eps}{2}\right)\geq
1+r-\mathbf{P}_x(\Omega_\eps).\end{equation}
 From (\ref{asintotica}), letting $t\rightarrow +\infty$, we get
\[V(x)\geq
\int_{B(x)}\int_0^{+\infty}\!\!l(\ox_s) ds  \geq
\int_{\Omega_\eps} \int_0^{+\infty}\!\!l(\ox_s)ds\geq
\int_{\Omega_\eps} \sum_{n} \int_{t_n} ^{t_{n}+h}\!\! l(\ox_s)
ds\geq
\]\[\geq \int_{\Omega_\eps} \sum_{n}h \inf_{[t_n(\omega),t_n(\omega)+h]}
l(\ox_t(\omega))\geq  h \sum_{n}\int_{\Omega_\eps}
\inf_{[t_n(\omega),t_n(\omega)+h]} l(\ox_t(\omega))\geq  \]
\[\geq  h \frac{\eps}{2} \sum_{n}\mathbf{P}
\left[\left(\sup_{t_n\leq s\leq t_n+h}
|l(\ox_{s})-l(\ox_{t_n})|\leq\frac{\eps}{2}\right)\cap\Omega_\eps\right]\geq
h \sum_{n}\frac{\eps}{2}r=+\infty\] where the last inequalities
are obtained using the strong Markov property of the process
$\ox_{\tau_K(.)}$.  This  gives a contradiction: then, for every
$\eps>0$, $\mathbf{P}(\Omega_\eps)=0$.
 We have proved that $l(\ox_t)\rightarrow 0$ as
$t\rightarrow +\infty$ for almost all $\omega\in B(x)$, now the
positive definiteness of $l$ implies that
\[\mathbf{P}_x\left\{\lim_{t\rightarrow+\infty}|\ox_t|=0\right\}\geq
\mathbf{P}_x(B_x)\geq 1-\frac{V(x)}{\min_{|y|\geq k}V(y)}.\]

We prove now the statement (ii). If $\mathcal{O}$ coincides with
the whole space, arguing  as above, we get that for every $k>0$
and $x\in B_k$ there exists a strict control $\alpha^k$ such that
the corresponding trajectory $X^k$ verifies
\[\mathbf{P}_x\left(\limsup_{t\rightarrow+\infty}l(X_t^k)>0\right)\leq
V(x)/\min_{|y|\geq k}V(y).\] Using the properness of the function
$V$, by letting $k\to +\infty$, we get that for every
$x\in\mathcal{O}$
\begin{equation}\label{infi}
\inf_{\alpha}\mathbf{P}_x\left(\limsup_{t\rightarrow+\infty}l(X_t^\alpha)>0\right)=0
\end{equation}
which gives, by the positive definiteness of the function $l$, the
asymptotic stabilizability in the large.
\end{proof}
\begin{rem}[Uniform asymptotic stabilizability in
probability]\upshape The existence of a Lyapunov function implies
a stronger asymptotic stability of the system, which we call
uniform asymptotic stabilizability. Moreover we will show in a
forthcoming paper that the uniform asymptotic stabilizability can
be completely characterized in terms of strict Lyapunov functions.

\noindent The system (CSDE) is {\em uniformly asymptotically
stabilizable in probability} in $\mathcal{O}$  if for every $x\in
\mathcal{O}$ there exists $\oa\in\mathcal{A}_x$ such that for
every $k>0$
\[\lim_{x\rightarrow 0}\mathbf{P}\left(\sup_{t\geq 0}|\ox_t|\geq
k\right)=0,\]\[\sup_{x\in\mathcal{O}}
T_x^{\oa}(\mathcal{O}\setminus B_k)<+\infty,\] where
$T_x^{\oa}(\mathcal{O}\setminus B_k)$ is the expected time spent
by the trajectory $\ox$ in the set $\mathcal{O}\setminus B_k$.

The fact that the existence of a Lyapunov function implies the
uniform asymptotic stabilizability follows very easily from the
representation formula for the function $V$ and the positive
definiteness of the function $l$:
\[ V(x)\geq\mathbf{E}_x\int_0^{\ot}\!\!l(\ox_s)ds
\geq T_x^{\oa}(\mathcal{O}\setminus B_r)\inf_{y\in
\mathcal{O}\setminus B_r}l(y),\] which implies
\[\sup_{x\in\mathcal{O}}T_x^{\oa}(\mathcal{O}\setminus B_r)\leq
\frac{\sup_{x\in\mathcal{O}}V(x)}{\inf_{y\in \mathcal{O}\setminus
B_r}l(y)}=\frac{\|V\|_\infty} {\inf_{y\in \mathcal{O}\setminus
B_r}l(y)}. \]

The proof of the fact that uniform asymptotic stability implies
asymptotic stability (in particular that for every initial data
there exists a control driving asymptotically the trajectory to
the origin almost surely) is an argument based on continuity
properties of trajectories of (CSDE) of the type
(\ref{continuitytrajectory}) we proved in Theorem \ref{liap2}.
\end{rem}

\section{Attractors}
Next we  extend the results in section 4 to study the
stabilizability of general closed sets $M\subseteq \R^N$. We
denote by $d(x,M)$  the distance between a point $x\in\R^N$ and
the set $M$.

We recall that a closed set $M$ is {\it viable} with respect to a
stochastic controlled dynamical system if for every $x\in M$ there
exists an admissible  control such that the corresponding
trajectory remains almost surely in $M$.

\begin{defi}[Stabilizability in probability at $M$]
A closed set $M$ is \emph{stabilizable in probability} for
$(CSDE)$ if for every $k>0$ there exists $\delta>0$ such that, for
every $x$  at distance less than $\delta$ from $M$, there exists
an admissible control $\oa$ such that the corresponding trajectory
$\ox_.$ verifies
\[\mathbf{P}_x\left(\sup_{t\geq 0}d(\ox_t,M)\geq k\right)\leq
\eps.
\]
\end{defi}

\begin{rem}\upshape
We observe that if $M$ is stabilizable in probability according to
the previous definition, then in particular it is viable. In fact
for every $\eps>0$ fixed, the definition gives that, for $x\in M$,
$\inf_\alpha\mathbf{E}_x\int_0^{+\infty }e^{-\lambda
t}k_\eps(X_t)dt=0$ for any $\lambda>0$ and any smooth function
$k_\eps$ which is nonnegative, bounded and null on the points at
distance less than $\eps$ from $M$. By  Theorem
\ref{optimalcontrol}, the infimum is attained, therefore for every
$\eps>0$ there is a control $\oa\in\mathcal{A}_x$ whose
corresponding  trajectory  stays almost surely at distance less
than $\eps$ from $M$: in particular, for every $\lambda>0$,
$\inf_\alpha\mathbf{E}_x\int_0^{+\infty} e^{-\lambda t}|X_t|dt=0$.
Therefore, again by Theorem \ref{optimalcontrol}, there exists,
for every $x\in M$, a minimizing control whose corresponding
trajectory stays in $M$ almost surely for all $t\geq 0$.

A geometric characterizations of viability of closed sets with
respect to a stochastic differential controlled equation has been
given in ~\cite{bj} (see also references therein). According to
this characterization, we note that the fact that the set $M$ is
stabilizable in probability implies that the diffusion has to
degenerate on its boundary: for every $x\in
\partial M$ there exists $\alpha\in A$  such that
$\sigma(x,\alpha)\cdot p=0$ for every $p$  generalized normal
vector to $M$ at $x$.
\end{rem}
We introduce the notion of controlled attractiveness: it
coincides, when the system is uncontrolled, with the standard
notion of pathwise forward attractiveness (see \cite{has}).
\begin{defi}[Controlled attractor]
The set $M$ is a {\em controlled attractor} for the system (CSDE)
in the open set $\mathcal{O}\subseteq \R^N$ if for every initial
data $x\in \mathcal{O}$ then \[\inf_\alpha
\mathbf{P}_x\left(\limsup_{t\rightarrow +\infty}d(X_t^\alpha, M)>0
\right)=0.\] This means that for every $\eps>0$ there exists
$\oa\in\mathcal{A}_x$ such that the corresponding  trajectory
approaches asymptotically the set $M$ with probability at least
$1-\eps$.\\
The set $\mathcal{O}$ is called {\em domain of attraction} for
$M$: if it coincides with $\R^N$ the set $M$ is a {\em global}
attractor.
\end{defi}

\begin{rem}\upshape We consider a
function $V:\R^N\rightarrow\R$ which satisfies the conditions  in
the Definition \ref{gsliap} of strict global Lyapunov function
with the only difference that the function $l$ is assumed only
nonnegative definite.  The proof of Theorem \ref{liap2} can be
repeated in this case: we obtain that  for every $x\in\R^N$
\begin{equation}\label{lasa}
\inf_{\alpha}\mathbf{P}_x\left(\limsup_{t\rightarrow+\infty}l(X_t^\alpha)>0\right)=0.\end{equation}
We introduce the set ${\cal L}:=\{y\ | \ l(y)=0\}$. From
(\ref{lasa}) we get that for every $x\in\R^N$
\[\inf_{\alpha}\mathbf{P}_x\left(\limsup_{t\rightarrow+\infty}d(X_t^\alpha,\mathcal{L})>0\right)=0,\]
which means that ${\cal L}$ is a controlled global attractor for
the system. For uncontrolled diffusion processes results of this
kind can be found in \cite{mao} and \cite{deng}. The earlier paper
of Kushner \cite{ku2} studies also a stochastic version of the La
Salle invariance principle, namely, that the omega limit set of
the process is an invariant subset of ${\cal L}$, in a suitable
sense.
\end{rem}

We can generalize the notion of control Lyapunov function in order
to study the attractiveness and the stabilizability of a set $M$.
\begin{defi}[control $M$-Lyapunov function]
Let $M\subseteq \R^N$ be a closed set  and
$\mathcal{O}\subseteq\R^N$ an open set containing $M$. A function
$V:\mathcal{O}\rightarrow [0,+\infty)$ is a \emph{control
M-Lyapunov function} for $(CSDE)$ if it satisfies

\noindent (i) it is lower semicontinuous and continuous at every
$x\in\partial M$;

\noindent(ii) it is {\em M-positive definite}, i.e., $V(x)>0$ for
$x\not\in M$ and $V(x)=0$ for $x\in M$;

\noindent(ii) $V$ is {\em M-proper}, i.e., its level sets  $\{x\
|\ V(x)\leq\mu\}$ are bounded;

\noindent(iii) it is a viscosity supersolution of the equation
\[
\max_{\alpha\in A}\left\{-DV(x)\cdot
f(x,\alpha)-trace\left[a(x,\alpha) D^2V(x)\right]\right\}\geq l(x)
\ \ \ \ x\in\mathcal{O}.\] If $l(x)\geq 0$ then $V$ is a control
$M$-Lyapunov function, if $l(x)$ is a $M$ positive definite,
Lipschitz continuous bounded function   then $V$ is a strict
control Lyapunov function.
\end{defi}
We can therefore prove for the case of a set $M$ very similar
results as for the case of an equilibrium point.
\begin{teo}
If the system $(CSDE)$ admits a control $M$-Lyapunov function $V$
then the system is $M$ stabilizable; if moreover the function $V$
is a strict control $M$-Lyapunov function then the set $M$ is a
controlled attractor for the system with domain of attraction
equal to $\mathcal{O}$.\end{teo}
\begin{proof}
In order to prove this result, one can repeat the proofs given in
Theorems \ref{liap1} and \ref{liap2}   since the function $V$
satisfies a superoptimality principle and the level sets of $V$
are a local basis of neighborhoods of $M$.
\end{proof}
\section{Examples}
In this section we present some very simple examples illustrating
the theory.

The first example is about a {\em stochastic perturbations of
stabilizable systems}. We apply the Lyapunov Theorems to show that
an asymptotically controllable deterministic dynamical system
continues to be stabilizable or asymptotically stabilizable in
probability if we perturb it with a white noise of  intensity
small enough. The idea to prove it relies on the fact that, if the
stochastic perturbation is small enough, then a Lyapunov function
for the deterministic systems remains a Lyapunov function also for
the stochastic one.

\begin{ex}\upshape We consider a deterministic controlled system in $\R^N$
\begin{equation}
\label{detsyst2} \dot X_t=f(X_t
,\alpha_t) 
\end{equation}
where $f(x,a)$ is a Lipschitz continuous, locally bounded function
in $x$ uniformly with respect to $a$ and the control $\alpha$ is a
measurable function taking values in a compact space $A$. We
assume  that  the system is globally asymptotically (open loop)
stabilizable at the origin, i.e., asymptotically controllable in
the  terminology of deterministic systems \cite{son, sosu}. By the
converse Lyapunov Theorem \cite{son0, sosu}, there exists a
continuous control Lyapunov function for the system, i.e., for
some positive definite continuous function $L$, there exists a
proper, positive definite function $V$ satisfying in $\R^N$
\begin{equation}
\label{ldet2} \max_{\alpha\in A}\left\{ -f(x,\alpha)\cdot
DV\right\} \geq L(x)
\end{equation}
in the viscosity sense.   Moreover we can choose the function $V$
to be semiconcave away from the origin as proved by Rifford in
\cite{r}. This means  that for every $\delta>0$  there exists a
semiconcavity constant $C_\delta>0$ such that the function
\[V(x)- \frac{C_\delta}{2} |x|^2\] is concave in $\R^N\setminus B_\delta$.
The semiconcavity constant $C_\delta$ is an upper bound on the
second derivatives of the function (to be intended in the sense of
distributions). In particular, by the definition of semiconcavity,
we get that if $|x|>\delta$ and $(p,X)\in \mathcal{J}^{2,-}V(x)$
then
\begin{equation}\label{semiconca} C_\delta \mathbf{I}_N-X \geq
0.\end{equation}

We  study under which conditions the system continues to be
asymptotically or Lyapunov stabilizable if we perturb it with  a
white noise. We consider the perturbed system
\[dX_t=f(X_t,\alpha_t)dt+\sigma(X_t,\alpha_t)dB_t\]
where $(B_t)_t$ is a $M$-dimensional white noise and the function
$\sigma(x,a)$ is Lipschitz continuous in $x$ uniformly  with
respect to $a$ and takes values in the space of the $N\times M$
dimensional matrices with entries in $\R$.

By the semiconcavity inequality (\ref{semiconca}) and by
(\ref{ldet2}), we get for every $|x|>\delta$, $(p,X)\in
J^{2,-}V(x)$
\[\max_{\alpha\in A}\left\{-f(x,\alpha)\cdot p-trace\ a(x,\alpha)X
\right\}\geq\max_{\alpha\in A}\left\{-f(x,\alpha)\cdot
p\right\}-\max_{\alpha\in A} \left\{trace\
a(x,\alpha)X\right\}\geq
\]\[\geq L(x)-C_\delta\max_{\alpha\in A} \left\{trace\
a(x,\alpha)\right\}.\] Therefore, if the diffusion $\sigma$
satisfies a small intensity condition
\[trace\ a(x,\alpha)\leq \frac{L(x)}{C_\delta}\quad \forall \alpha\in A \ \forall \ |x|>\delta,\]
we can conclude  that  the function $V$ is a control Lyapunov
function for the stochastic system and then, according to Theorem
\ref{liap1}, the system is both Lyapunov and Lagrange stabilizable
in probability.

If moreover   for every $\delta>0$
\[trace\ a(x,\alpha)<\frac{L(x)}{C_\delta}
\quad \forall \alpha\in A \ \forall\  |x|>\delta,\] it is possible
to construct a positive definite, Lipschitz continuous function
$l$ such that $V$ is a viscosity supersolution of
\[\max_{a\in A}\left\{-f(x,a)\cdot
DV(x)-trace\ a(x,\alpha)D^2V(x)\right\} \geq l(x)\] and then, by
Theorem \ref{liap2}, the system is asymptotically stabilizable in
the large at the equilibrium.

A similar result can be obtained in the case of local
asymptotically controllable systems.
\end{ex}
In the next example we  give conditions on a radial function to be
a Lyapunov function for the stability in probability.
\begin{ex}\upshape In this example we consider
as candidate Lyapunov function for the general controlled system
$(CSDE)$ the function $V(x)=|x|^\gamma$ for some positive
$\gamma>0$ and study under which conditions the system is
stabilizable.

We compute \[DV(x)=\gamma |x|^{\gamma-2}x \quad  D^2V(x)=\gamma
|x|^{\gamma-2} \mathbf{I}+\gamma(\gamma-2) |x|^{\gamma-4}(x\cdot
x^T).\] Therefore \[\max_{a\in A}\left\{-f(x,a)\cdot DV(x)-trace\
a(x,\alpha)D^2V(x)\right\}=\]\[=\gamma |x|^{\gamma-2}\max_{a\in
A}\left\{-f(x,a)\cdot x -trace\ a(x,\alpha)-(\gamma-2)
|x|^{-2}trace\ a(x,\alpha)(x\cdot x^T))\right\}=\]\[= \gamma
|x|^{\gamma-2}\max_{a\in A}\left\{-f(x,a)\cdot x -trace\
a(x,\alpha)-(\gamma-2) |x|^{-2} |\sigma(x,a)^T\cdot
x|^2\right\}.\] If $\gamma\leq 2$, this gives that the $V$ is a
Lyapunov function for the system if for every $x$ there exists
$\alpha\in A$ such that $f(x,\alpha)\cdot x+trace\ a(x,\alpha)\leq
0$. We can observe that, since $trace\, a(x,\alpha)\geq 0$ for
every $\alpha$, the radial component of the drift $f$ must be
everywhere nonpositive, for some $\alpha\in A$. In particular it
must be negative to compensate the destabilizing role of the
diffusion, whenever $trace\, a(x,\alpha)$ is nonnull.
\end{ex}

\end{document}